\documentclass[11pt]{amsart}

\usepackage[utf8]{inputenc}
\usepackage[T1]{fontenc}
\usepackage{lmodern}
\usepackage{amsmath,amssymb,amsthm}
\usepackage[margin=1.15in]{geometry}
\usepackage{microtype}
\usepackage{csquotes}
\usepackage{xcolor}
\usepackage[colorlinks=true,linkcolor=blue!50!black,citecolor=blue!50!black,urlcolor=blue!50!black]{hyperref}

\newtheorem*{theorem*}{Theorem}

\title[Truth, Beauty, and the Machine]{Truth, Beauty, and the Machine: Why Mathematicians Will Not Become Irrelevant?}

\author{Ashish K. Srivastava}
\address{Department of Mathematics and Statistics, Saint Louis University, St.~Louis, MO, USA}
\email{ashish.srivastava@slu.edu}

\begin{document}

\begin{abstract}
Machines can now solve the most difficult open problems, and many mathematicians, especially young ones, are asking whether the discipline still has a future. Mathematics has never been simply the production of true statements. It is a human search for truth that is \emph{understood}, and for beauty that tells us we have understood it. A machine that settles a famous conjecture enlarges the stock of known truths, but it does not remove the need for someone to see why the statement is true, to judge whether it matters, and to fold it into a structure that others can live in. These are not leftover tasks waiting for better software. They are what mathematics is.
\end{abstract}

\maketitle

\section{A question that everyone is asking}

\noindent In the past two months, mathematics community has been rocked by AI solving two of the seven Millenium Prize Problems selected by the Clay Mathematics Institute: the Riemann Hypothesis, and the Navier-Stokes Problem. These problems have been considered the most difficult open problems in the field of mathematics. In addition to these two giant leaps, several hundred other considerably difficult problems have been solved in just past two-three months including Jacobian Conjecture, Koethe conjecture and pure-semisimplicity conjecture. These developments have been unsettling for senior researchers and young students alike. In the past few weeks, undergraduate and graduate students in my abstract algebra class have been asking me difficult questions such as: do you still feel motivated to do research in mathematics? Do you still think we have a future in mathematics? This essay is my attempt to answer those difficult questions. In 1888 Richard Dedekind published a small book with the title: \emph{Was sind und was sollen die Zahlen?}, roughly translated in English as ``What are numbers, and what should they be?'' \cite{Dedekind}. In this book he characterized the natural numbers by axioms rather than by calculation, and so helped start the move toward abstraction that shaped twentieth-century mathematics. The move was resisted, sometimes bitterly. To many working mathematicians of the time, it seemed to drain the subject of its concrete content but slowly the mathematics community got used to it and accepted it as an advancement. I am of the view that mathematics will survive this transition, as it survived the Dedekind transition, and in this essay I want to press a related point more directly, perhaps more stubbornly. Mathematics is a pursuit of truth and of beauty, and pursuits of that kind do not become irrelevant when a new instrument arrives. The telescope did not make astronomers obsolete. It gave them more sky.

Today no one can doubt the ability of AI in solving problems. I expect that they will solve many more in coming days, and some of them will be problems we have cared about for generations. My claim is that solving problems was never the whole of the matter, and that the parts which remain are the parts that give mathematics its meaning.

\section{What a proof is for}

\noindent Mathematics is the one human activity that claims to reach
certainty, and the proof is the means by which that certainty is reached. One
might therefore conclude that once a machine can produce proofs, and especially
once those proofs are checked by a formal verifier, the essential work is done.
The truth has been established. What more could anyone want?

The history of the subject answers this question clearly: we want to
\emph{understand}. William Thurston made this point in his essay \emph{On proof
and progress in mathematics} \cite{Thurston}. He argued that what mathematicians
produce is not a list of certified theorems but human understanding, and that
the value of a proof lies largely in what it allows others to see. A proof that
nobody can follow certifies a fact, but it does not yet hand that fact over to
the community as knowledge in the full sense.

We already have long experience with truths that arrived ahead of
understanding. The Four Color Theorem was proved in 1976 with essential help
from a computer \cite{AppelHaken}. The proof was correct, as later
re-examinations and an eventual formal verification have confirmed, and yet
the mathematical community did not regard the matter as finished. People
continued, and continue today, to look for a proof that explains \emph{why}
four colors suffice. Thomas Hales's proof of the Kepler conjecture was
similarly so dependent on computation that referees could not fully certify it;
the response was the Flyspeck project, a complete formal proof completed by a
large collaboration \cite{Hales}. In both cases the arrival of machine-assisted
truth did not close the subject. It opened new questions, about
explanation, about verification, and about what it means to know.

There is a useful distinction here between \emph{knowing that} and
\emph{knowing why}. A machine-generated proof of a deep conjecture gives us the
first. The second requires that someone, somewhere, reorganize their
understanding so that the result becomes natural, even inevitable. Often this
requires new concepts. The proof of Fermat's Last Theorem matters to number
theory less as a resolution of Fermat's marginal note than as a demonstration of
the power of modularity, a body of ideas that now organizes a great deal of the
subject. If a machine had produced a valid but conceptually opaque proof of
Fermat in 1990, the community would have learned that the theorem was true and
would still have had to discover everything that actually made the proof
important.

Truth in mathematics, then, has two layers. There is truth as a verdict:
this statement holds. And there is truth as insight: this is how the world of
mathematical objects is arranged, and here is why it could not be otherwise.
Machines are becoming very good at the first. The second is not a product that
can be delivered. It is a state that a mind has to reach.

\section{Beauty is not decoration}

\noindent G.~H. Hardy wrote that there is \enquote{no permanent place in the world for
ugly mathematics} \cite{Hardy}. The remark is often quoted as aesthetic
snobbery, but I think it is better read as a claim about knowledge. Beauty in
mathematics is not ornament applied after the fact. It is one of the most
reliable signals, a kind of litmus test we have that we have found the right way to see something. 

Consider what mathematicians mean when they call a proof beautiful or in the words of Paul Erd\"{o}s, a proof straight from book. Usually
they mean that it is economical, that it reveals a hidden connection, and that
it makes the result seem inevitable once the key idea is in view. Euclid's proof
that there are infinitely many primes is beautiful because one simple
construction, multiplying the known primes and adding one, does all the work.
Such a proof is short not because it skips steps but because it has found the
joint at which the problem divides.

Let me give an example from my own research area of noncommutative algebra, because it shows how
beauty and truth are bound together in practice. Most rings we meet in a first
course have the \emph{invariant basis number} property: if $R^m \cong R^n$ as
$R$-modules, then $m = n$. This is the property that makes rank of a free module well-defined. In 1962 William Leavitt showed that this can fail in an extreme and
precise way \cite{Leavitt}. For each $n \ge 2$ there is an algebra
$L = L_K(1,n)$ over a field $K$ for which
\[
  L \;\cong\; L^{n}
\]
as left $L$-modules. The algebra is generated by elements
$x_1,\dots,x_n,y_1,\dots,y_n$ subject to the relations
\[
  y_i x_j = \delta_{ij}, \qquad \sum_{i=1}^{n} x_i y_i = 1,
\]
which say, in effect, that a row and a column of length $n$ are inverse to each
other. Leavitt's algebras might have remained a curiosity. Instead, about four
decades later, it was recognized that they are the simplest members of a large
family: to every directed graph $E$ one can attach a \emph{Leavitt path algebra}
$L_K(E)$, and $L_K(1,n)$ is the algebra of the graph with one vertex and $n$
loops \cite{AbramsAranda,AraMorenoPardo}. Ring-theoretic properties that seemed
hard to control, among them simplicity, pure infiniteness, GK-dimension, and the structure of
the monoid of finitely generated projective modules, can then be read off from
the combinatorics of the graph.

This is what mathematical beauty looks like from the inside. A strange
phenomenon becomes an instance of a pattern; the pattern turns out to be
governed by a simple picture; and the picture links algebra to
operator algebras, to symbolic dynamics, and to $K$-theory. None of this is
needed to \emph{certify} Leavitt's original theorem. All of it is needed to
\emph{understand} it. The beauty was not a reward handed out after the
truth had been established. It was the route by which the truth became
comprehensible, and it is what made the subject grow.

Keats's line, \enquote{Beauty is truth, truth beauty}, is too simple as
philosophy. As a description of how mathematicians actually work, it is close
to right. Henri Poincar\'e argued that the mathematician's aesthetic sense
acts as a filter, singling out from the vast number of possible combinations
the few that are fruitful \cite{Poincare}. If that is so, then aesthetic
judgment is not a luxury that machines will make unnecessary. It is the
faculty that decides which of the machine's outputs are worth our
attention.

\section{The chess player and the mathematician}

\noindent An analogy is often offered at this point, and it is worth taking seriously,
along with its limits. Chess engines have exceeded the best human players for
more than a quarter of a century. Yet chess is not dead. More people play it now
than ever. The engines have become teachers, sparring partners, and sources of
new ideas, and human games remain meaningful because they are human.

The analogy shows that being surpassed by machines at a task does not abolish
the value of a human practice. But mathematics is not a game, and I want to
resist the suggestion that it will survive only as a kind of high-level hobby,
admirable but inessential. The difference is this. A chess game has no
content beyond itself; its value is in the playing. A mathematical theorem has
content that reaches beyond the person who proved it, and its value depends on
that content being absorbed into a shared understanding. Mathematics is
therefore both a practice, like chess, and a body of knowledge, like
science. Machines change the practice. They cannot make the knowledge
unnecessary, because knowledge in the full sense is not something that exists
on a disk. It exists when someone understands.

Edsger Dijkstra once suggested that asking whether a computer can think is
about as interesting as asking whether a submarine can swim, a remark I found as a footnote in Avigad's essay \cite{Avigad}. The submarine crosses the ocean; it does not swim,
and nobody is troubled by this. Similarly, I do not need to settle whether a
machine ``understands'' a proof to observe that \emph{my} understanding of it
is not supplied by the machine's having produced it. Understanding is not
transferable by file. It has to be rebuilt in each mind that comes to possess
it, and that rebuilding is precisely what mathematical education, mathematical
exposition, and mathematical research consist of.

\section{Why mathematics will not become irrelevant}

\noindent Every proof answers a question, and somebody had to ask it. Choosing what
to ask, what to define, and what to regard as important is not a
computational task with a correct output. It is a judgment of value. The
definition of a group, of a topological space, of a sheaf, of a Leavitt path
algebra: each was a decision that a certain pattern deserved a name and a
theory. Machines can propose definitions, and they may become very good
collaborators in doing so. But the decision about which definitions to live
with, to teach, and to build on is a decision a community of people makes about
how it wants to understand the world.

A result that nobody understands is, practically speaking, a sealed
box. To use a theorem, whether in further mathematics or in physics or
engineering, one has to know what its hypotheses mean, when they hold, and
what the conclusion actually says. This interpretive work does not shrink as
machines prove more. It grows. The more truths are produced, the greater the
need for people who can read them, connect them, and judge their significance.

It is worth remembering that the systems now proving theorems are themselves
built from linear algebra, probability, optimization, and the theory of
computation. Questions about why these systems work, when they fail, and how
far their outputs can be trusted are mathematical questions, and many of them are
open. Eugene Wigner wrote of the \enquote{unreasonable effectiveness} of
mathematics in the natural sciences \cite{Wigner}. We may now add its effectiveness in explaining
our own inventions. A world increasingly run by algorithms needs more
mathematical understanding, not less.

Finally, mathematics is how human beings learn to reason with care: to define
terms, to separate assumption from conclusion, to tell the difference between
a convincing argument and a correct one. That last distinction matters more
than ever in an age of fluent machines that can generate arguments that
merely sound right. A society that has handed its proofs to machines still
needs citizens who know what a proof is.

\section{What should change, and what should not}

\noindent The advancement in AI we have seen in the past few months has made one thing clear that certain kinds of technical labor that once occupied years of a career will be done quickly. Some problems we hoped to solve ourselves will be solved by machines first. Young mathematicians in particular deserve to hear this plainly, and also to hear the rest of the story: that the core of what they love, the struggle toward understanding and the moment when a structure suddenly becomes clear, was never about being the first to finish. A student who works through a proof that was found long ago is not doing something worthless because the theorem is already known. They are doing the thing that makes the theorem known \emph{to them}, which is the only way any theorem has ever been known to anyone.

What should not change are the values that make mathematics what it is: honesty about what has and has not been proved, care with definitions, the patience to seek the right idea rather than any idea that works, and
respect for beauty as a guide to truth. These are not old-fashioned sentiments to be traded for efficiency. They are exactly the standards by which we will have to judge the flood of results that machines are about to produce.

Mathematics is the human pursuit of truth that is understood and of beauty that signals understanding. Machines will become extraordinary instruments in that pursuit, and they may well reach some of its destinations before we do. But a destination reached is not the same as a place understood, and mathematics has always been about the understanding.

So long as there are people who want to know not only \emph{that} something is true but \emph{why}, and who feel the particular joy of seeing a hard thing become simple, mathematics will not be irrelevant. It will be what it has
always been: one of the deepest ways in which human beings come to know the world and themselves. I would like to end by quoting Flexner \cite{Flexner} ``We make ourselves no promises, but we cherish the hope that the unobstructed pursuit of useless knowledge will prove to have consequences in the future as in the past".


\begin{thebibliography}{99}

\bibitem{AbramsAranda}
G.~Abrams and G.~Aranda Pino,
\emph{The Leavitt path algebra of a graph},
J. Algebra \textbf{293} (2005), 319--334.

\bibitem{AppelHaken}
K.~Appel and W.~Haken,
\emph{Every planar map is four colorable. Part~I: Discharging},
Illinois J. Math. \textbf{21} (1977), 429--490.

\bibitem{AraMorenoPardo}
P.~Ara, M.~A. Moreno, and E.~Pardo,
\emph{Nonstable $K$-theory for graph algebras},
Algebr. Represent. Theory \textbf{10} (2007), 157--178.

\bibitem{Avigad}
J.~Avigad,
\emph{What is mathematics now, and what should it be?},
preprint (2026), \href{https://arxiv.org/abs/2608.23218}{arXiv:2608.23218}.

\bibitem{Dedekind}
R.~Dedekind,
\emph{Was sind und was sollen die Zahlen?},
Vieweg, Braunschweig, 1888.

\bibitem{Flexner}
A. Flexner, \emph{The usefulness of useless knowledge}, Harper's Magazine, 179 (1939).      

\bibitem{Hales}
T.~Hales et al.,
\emph{A formal proof of the Kepler conjecture},
Forum Math. Pi \textbf{5} (2017), e2.

\bibitem{Hardy}
G.~H. Hardy,
\emph{A Mathematician's Apology},
Cambridge University Press, 1940.

\bibitem{Leavitt}
W.~G. Leavitt,
\emph{The module type of a ring},
Trans. Amer. Math. Soc. \textbf{103} (1962), 113--130.

\bibitem{Poincare}
H.~Poincar\'e,
\emph{Science et m\'ethode},
Flammarion, Paris, 1908.

\bibitem{Thurston}
W.~P. Thurston,
\emph{On proof and progress in mathematics},
Bull. Amer. Math. Soc. (N.S.) \textbf{30} (1994), 161--177.

\bibitem{Wigner}
E.~P. Wigner,
\emph{The unreasonable effectiveness of mathematics in the natural sciences},
Comm. Pure Appl. Math. \textbf{13} (1960), 1--14.


\end{thebibliography}
\end{document}